\documentclass[10pt,reqno]{amsart}
\usepackage{amssymb}
\usepackage{amsmath}
\usepackage{amscd}
\begin{document}
\setlength{\oddsidemargin}{0cm} \setlength{\evensidemargin}{0cm}

\theoremstyle{plain}
\newtheorem{theorem}{Theorem}[section]
\newtheorem{proposition}[theorem]{Proposition}
\newtheorem{lemma}[theorem]{Lemma}
\newtheorem{corollary}[theorem]{Corollary}
\newtheorem{conj}[theorem]{Conjecture}

\theoremstyle{definition}
\newtheorem{definition}[theorem]{Definition}
\newtheorem{exam}[theorem]{Example}
\newtheorem{remark}[theorem]{Remark}

\numberwithin{equation}{section}

\title[]{Non-Riemannian Einstein-Randers metrics on $E_6/A_4$ and $E_6/A_1$}

\author{Xiaosheng Li}
\address{School of Mathematics and Statistics, Xinyang Normal University, Xinyang 464000, P.R. China} \email{xiaosheng0526@126.com}

\author{Chao Chen}
\address{School of Mathematics and Statistics, Xinyang Normal University, Xinyang 464000, P.R. China} \email{}

\author{Zhiqi Chen}
\address{School of Mathematical Sciences and LPMC, Nankai University, Tianjin 300071, P.R. China} \email{chenzhiqi@nankai.edu.cn}

\author{Yuwang Hu}
\address{School of Mathematics and Statistics, Xinyang Normal University, Xinyang 464000, P.R. China} \email{hywzrn@163.com}



\begin{abstract}
In this paper, we first prove that homogeneous spaces $E_6/A_4$ and $E_6/A_1$ admit Einstein metrics which are $Ad(T\times A_1\times A_4)$-invariant, and then show that they admit Non-Riemannian Einstein-Randers metrics.
\end{abstract}

\maketitle


\setcounter{section}{0}
\section{Introduction}
Randers metrics were introduced by Randers in the context of general
relativity, and named after him by Ingarden. It shows the importance
of Randers metrics in physics. Moreover, Randers metrics are useful
in other fields; see Ingarden's account on \cite{AIM,In1} for their
application in the study on the Lagrangian of relativistic electrons.

Just as in the Riemannian case, it is a fundamental problem to
classify homogeneous Einstein-Randers spaces. In particular, it is
very important to know if a homogeneous manifold admits invariant Einstein-Randers metrics. There are a lot of studies on Einstein-Randers metrics on homogeneous manifolds, see \cite{CDL1,CDL2,KC1,LD1,NT1,WD1,WD2,WD3,WHD1}.

In this paper, we will discuss Einstein metrics and Einstein-Randers metrics on homogeneous spaces $E_6/A_4$ and $E_6/A_1$. In Section 2, we prove that there are only four Einstein metrics on $E_6/A_4$ and two Einstein metrics on $E_6/A_1$ which are $Ad(T\times A_1\times A_4)$-invariant. Furthermore, in Section 3, we prove that there are at least four and two families of $E_6$-invariant non-Riemannian Einstein-Randers metrics on $E_6/A_4$ and $E_6/A_1$ respectively.

\section{Einstein metrics}
Consider the symmetric spaces $E_6/A_1\times A_5$ and $A_5/T\times A_4$. Let ${E_6}={A_1}\oplus{A_5}\oplus{\mathfrak m}_1$ and $A_5={\mathfrak h_0}\oplus A_4\oplus {\mathfrak m}_2$ be the corresponding decompositions of the Lie algebras. Here $\mathfrak h_0$ is the Lie algebra of $T$. Let $H$ be the Lie group $T\times A_4\times A_1$ with the Lie algebra ${\mathfrak h}_0\oplus A_4\oplus A_1$. Then we have the following decomposition of the Lie algebra $E_6$:
\begin{equation}\label{dec}
{E_6}={\mathfrak h}_0\oplus{A_4}\oplus{A_1}\oplus{\mathfrak m}_1\oplus{\mathfrak m}_2.
\end{equation}
Here the $Ad(H)$-modules ${\mathfrak m}_i,i=1,2$ are irreducible and
mutually non-equivalent and $\dim {\mathfrak
h}_0=1$. For the structure of such decomposition, see \cite{AMS1,DK1}. Clearly $\dim {A_4}=24$, $\dim {A_1}=3$, $\dim {\mathfrak m}_1=40$ and $\dim {\mathfrak m}_2=10$. The left-invariant metric on $E_6$ which
is $Ad(H)$-invariant must be of the form
\begin{equation}\label{metricII}
\langle\cdot,\cdot\rangle=u_0\cdot B|_{\mathfrak h_0}+u_1\cdot
B|_{A_4}+u_2\cdot B|_{A_1}+x_1\cdot
B|_{\mathfrak m_1}+x_2\cdot B|_{\mathfrak m_2}, \end{equation} where
$u_0,u_1,u_2,x_1,x_2\in {\mathbb R}^+$, and the space of left-invariant symmetric covariant 2-tensors on $E_6$ which are
$Ad(H)$-invariant is given by
\begin{equation}\label{non-naturalII}
v_0\cdot B|_{\mathfrak h_0}+ v_1\cdot B|_{A_4}+v_2\cdot
B|_{A_1}+v_3\cdot B|_{\mathfrak m_{1}}+v_4\cdot
B|_{\mathfrak m_{2}},
\end{equation}
where $v_0, v_1, v_2, v_3, v_4\in {\mathbb R}$. In particular, the
Ricci tensor $r$ of a left-invariant Riemannian metric $\langle\cdot,\cdot\rangle$ on $G$ is a left invariant symmetric covariant 2-tensor on
$G$ which is $Ad(H)$-invariant. Thus $r$ is of the form
(\ref{non-naturalII}). The authors give in \cite{AMS1} the formulae of the Ricci tensor corresponding to the metric~(\ref{metricII}) on $E_6$. In fact, the case for $E_6$ is just one case of what the authors study in \cite{AMS1}.

For this decomposition of $E_6$, the left-invariant metric on $E_6/A_4\times A_1$ which is $Ad(H)$-invariant must be of the form
\begin{equation}\label{metric2II}
\langle\cdot,\cdot\rangle=u_0\cdot B|_{\mathfrak h_0}+x_1\cdot
B|_{\mathfrak m_1}+x_2\cdot B|_{\mathfrak m_2},
\end{equation} where $u_0,x_1,x_2\in {\mathbb R}^+$. Based on the formulae given in \cite{AMS1}, the authors classify in \cite{CDL2} Einstein metrics on $E_6/A_4\times A_1$ which are $Ad(H)$-invariant.

The following is to discuss the left-invariant metrics on $E_6/A_4$ and $E_6/A_1$ which are $Ad(H)$-invariant corresponding to the decomposition~(\ref{dec}) of $E_6$.

\subsection{The case of $E_6/A_4$}
The left-invariant metric on $E_6/A_4$ which is $Ad(H)$-invariant must be of the form
\begin{equation}\label{met1}
\langle\cdot,\cdot\rangle=u_0\cdot B|_{\mathfrak h_0}+u_2\cdot B|_{A_1}+x_1\cdot
B|_{\mathfrak m_1}+x_2\cdot B|_{\mathfrak m_2}, \end{equation} where
$u_0,u_2,x_1,x_2\in {\mathbb R}^+$. Based on the formulae given in \cite{AMS1}, we have the components of the Ricci tensor $\widetilde{r}$ of the metric (\ref{met1}) on $E_6/A_4$:
\[ \left\{ \begin{aligned}
&
  \widetilde{r}_{{\mathfrak h}_0}=\frac{u_0}{8x_1^2}+\frac{u_0}{8x_2^2},\\
   &
  \widetilde{r}_{A_1}=\frac{1}{24u_2}+  \frac{5u_2}{24x_1^2}£¬ \\
  &
  \widetilde{r}_{{\mathfrak
  m}_1}=\frac{1}{2x_1}-\frac{x_2}{16x_1^2}-\frac{u_0}{160x_1^2}-\frac{u_2}{32x_1^2},\\
  &
  \widetilde{r}_{{\mathfrak
  m}_2}=\frac{1}{4x_2}+\frac{x_2}{8x_1^2}-  \frac{u_0}{40x_2^2}.
\end{aligned} \right. \]
Furthermore, the metric is Einstein if and only if there
exists a positive solution $\{u_0, u_2, x_1, x_2\}$ of the system of
equations
\begin{equation}\label{1}
\widetilde{r}_{{\mathfrak h}_0}=\widetilde{r}_{A_1}=\widetilde{r}_{{\mathfrak
  m}_1}=\widetilde{r}_{{\mathfrak   m}_2}.
\end{equation}
The following is to solve the equations by the theory of Gr$\rm{\ddot{o}}$bner basis. Putting $u_0=1$ and by $\widetilde{r}_{{\mathfrak h}_0}=\widetilde{r}_{A_1}, \widetilde{r}_{{\mathfrak h}_0}=\widetilde{r}_{{\mathfrak  m}_1}, \widetilde{r}_{{\mathfrak h}_0}=\widetilde{r}_{{\mathfrak   m}_2}$, we have
\[ \left\{ \begin{aligned}
  & f_1=-x_1^2x_2^2-5x_2^2u_2^2+3x_1^2u_2+3x_2^2u_2=0, \\
  & f_2=-80x_1x_2^2+10x_2^3+20x_1^2+21x_2^2+5x_2^2u_2=0, \\
  & f_3=-10x_1^2x_2-5x_2^3+6x_1^2+5x_2^2=0.
\end{aligned} \right. \]
Consider the polynomial ring $R={\mathbb Q}[z, x_1, x_2, u_2]$ and an ideal $I$ generated by $\{f_1, f_2, f_3,
zx_1x_2u_2-1\}$ to find non-zero solutions of $(\ref{1})$. Take a lexicographic order $>$ with
$z >u_2> x_1> x_2$ for a monomial ordering on $R$. By the help of computer, we have the polynomial of $x_2$ containing in the Gr$\rm{\ddot{o}}$bner basis of the ideal $I$:
\begin{eqnarray*}
  f(x_2)&=& 27263765625x_2^8-82709987500x_2^7+94104102500x_2^6-48116787500x_2^5 \\
         && +9948352750x_2^4-491681700x_2^3+74376100x_2^2-1183780x_2+142129.
\end{eqnarray*}
In the Gr$\rm{\ddot{o}}$bner basis of the ideal $I$, $x_1$ and $u_2$ can be written into polynomials of $x_2$. The equation $f(x_2)=0$ has four solutions:
$$x_2\approx 0.6513015810, \quad x_2\approx0.6770950751, \quad x_2\approx0.8288266917, \quad x_2\approx0.8641265950.$$
In fact, we have the following solutions of $(\ref{1})$:
\begin{eqnarray*}
 \{u_2\approx0.1141930615, & x_1\approx1.200678505, & x_2\approx0.6513015810\}, \\
 \{u_2\approx1.746579387, & x_1\approx0.9798479028, & x_2\approx0.6770950751\}, \\
 \{u_2\approx0.7564861893, & x_1\approx0.5068895851, & x_2\approx0.8288266917\}, \\
 \{u_2\approx0.0549236976, & x_1\approx0.4382514353, & x_2\approx0.8641265950\}.
\end{eqnarray*}
In summary, we have:
\begin{theorem}\label{theorem1}
Let the notations be as above. Then every left invariant metric on
$E_6/A_4$ which is $Ad(H)$-invariant is of the form $(\ref{met1})$. Up to scaling, there are four Einstein metrics on $E_6/A_4$ which are $Ad(H)$-invariant.
\end{theorem}

\subsection{The case of $E_6/A_1$}
The left-invariant metric on $E_6/A_1$ which is $Ad(H)$-invariant must be of the form
\begin{equation}\label{met2}
\langle\cdot,\cdot\rangle=u_0\cdot B|_{\mathfrak h_0}+u_1\cdot B|_{A_4}+x_1\cdot
B|_{\mathfrak m_1}+x_2\cdot B|_{\mathfrak m_2}, \end{equation} where
$u_0,u_1,x_1,x_2\in {\mathbb R}^+$. Based on the formulae given in \cite{AMS1}, we have the components of the Ricci tensor $\widetilde{r}$ of the metric (\ref{met2}) on $E_6/A_1$:
\[ \left\{ \begin{aligned}
&
  \widetilde{r}_{{\mathfrak h}_0}=\frac{u_0}{8x_1^2}+\frac{u_0}{8x_2^2},\\
   &
  \widetilde{r}_{A_4}=\frac{5}{48u_1}+  \frac{u_1}{8x_1^2}+\frac{u_1}{48x_2^2}£¬ \\
  &
  \widetilde{r}_{{\mathfrak
  m}_1}=\frac{1}{2x_1}-\frac{x_2}{16x_1^2}-\frac{u_0}{160x_1^2}-\frac{3u_1}{20x_1^2},\\
  &
  \widetilde{r}_{{\mathfrak
  m}_2}=\frac{1}{4x_2}+\frac{x_2}{8x_1^2}-  \frac{u_0}{40x_2^2}-\frac{u_1}{10x_2^2}.
\end{aligned} \right. \]
Furthermore, the metric is Einstein if and only if there
exists a positive solution $\{u_0, u_1, x_1, x_2\}$ of the system of
equations
\begin{equation}\label{2}
\widetilde{r}_{{\mathfrak h}_0}=\widetilde{r}_{A_4}=\widetilde{r}_{{\mathfrak
  m}_1}=\widetilde{r}_{{\mathfrak   m}_2}.
\end{equation}
Similar to the discussion on $E_6/A_4$. Putting $u_0=1$ and by $\widetilde{r}_{{\mathfrak h}_0}=\widetilde{r}_{A_4}, \widetilde{r}_{{\mathfrak h}_0}=\widetilde{r}_{{\mathfrak  m}_1}, \widetilde{r}_{{\mathfrak h}_0}=\widetilde{r}_{{\mathfrak   m}_2}$, we have
\[ \left\{ \begin{aligned}
  & f_1=-5x_1^2x_2^2-6x_2^2u_1^2-x_1^2u_1^2+6x_1^2u_1+6x_2^2u_1=0, \\
  & f_2=-80x_1x_2^2+10x_2^3+20x_1^2+21x_2^2+24x_2^2u_1=0, \\
  & f_3=-10x_1^2x_2-5x_2^3+6x_1^2+5x_2^2+4x_1^2u_1=0.
\end{aligned} \right. \]
Consider the polynomial ring $R={\mathbb Q}[z, x_1, x_2, u_1]$ and an ideal $I$ generated by $\{f_1, f_2, f_3,
zx_1x_2u_1-1\}$ to find non-zero solutions of $(\ref{2})$. Take a lexicographic order $>$ with
$z >u_1> x_1> x_2$ for a monomial ordering on $R$. By the help of computer, we have the polynomial of $x_2$ containing in the Gr$\rm{\ddot{o}}$bner basis of the ideal $I$:
\begin{eqnarray*}
  f(x_2)&=& 40733269776x_2^8-95717471616x_2^7+80248108328x_2^6-31589680504x_2^5 \\
        &&+7669961625x_2^4-1207801950x_2^3+120201725x_2^2-5089500x_2+422500.
\end{eqnarray*}
In the Gr$\rm{\ddot{o}}$bner basis of the ideal $I$, $x_1$ and $u_1$ can be written into polynomials of $x_2$. The equation $f(x_2)=0$ has two solutions:
$$x_2\approx0.8651778712,\quad x_2\approx 0.9203114422.$$
In fact, we have the following solutions of $(\ref{2})$:
\begin{eqnarray*}
 \{u_1\approx 0.1945580092,& x_1\approx 0.5189654864, & x_2\approx 0.8651778712\}, \\
 \{u_1\approx 0.7881276805,& x_1\approx 2.582407960, & x_2\approx 0.9203114422\}.
\end{eqnarray*}
In summary, we have:
\begin{theorem}\label{theorem2}
Let the notations be as above. Then every left invariant metric on
$E_6/A_1$ which is $Ad(H)$-invariant is of the form $(\ref{met2})$. Up to scaling, there are two Einstein metrics on $E_6/A_1$ which are $Ad(H)$-invariant.
\end{theorem}

\section{Einstein-Randers metrics}
A Randers metric $F$ on $M$ is built from a Riemannian metric and a
1-form, i.e.,
$$F=\alpha+\beta,$$ where $\alpha$ is a Riemannian metric and
$\beta$ is a 1-form whose length with respect to the Riemannian
metric $\alpha$ is less than 1 everywhere. Obviously, a Randers
metric is Riemannian if and only if it is reversible, i.e.,
$F(x,y)=F(x,-y)$ for any $x\in M$ and $y\in T_x(M)$. Sometimes it is
convenient to use the following presentation of a Randers metric in
\cite{BR1}, i.e.,
\begin{equation}
F(x,y)=\frac{\sqrt{[h(W,y)]^2+h(y,y)\lambda}}{\lambda}-\frac{h(W,y)}{\lambda}
\end{equation}
here $\lambda=1-h(W,W)>0$. The pair $(h,W)$ is called the navigation data of the
corresponding Randers metric $F$ .

The Ricci scalar $\mathfrak{Ric}(x,y)$ of a Finsler metric is
defined to be the sum of those $n-1$ flag curvatures $K(x,y,e_v)$,
where $\{e_v: v=1,2,\cdots,n-1\}$ is any collection of $n-1$
orthonormal transverse edges perpendicular to the flagpole, i.e.
\begin{equation}
 \mathfrak{Ric}(x,y)=\sum^{n-1}_{v=1}R_{vv}.
\end{equation}
The Ricci tensor is defined by
\begin{equation}
 Ric_{ij}=(\frac{1}{2}F^2\mathfrak{Ric})_{y^iy^j}.
\end{equation}
Obviously, the Ricci scalar depends on the position $x$ and the
flagpole $y$, but does not depend on the specific $n-1$ flags
with transverse edges orthogonal to $y$ (see \cite{BR1,BRS}). In the Riemannian case, it is a well known fact that the
Ricci scalar depends only on $x$. Thus it is
quite interesting to study a Finsler manifold whose Ricci scalar does
not depend on the flagpole $y$. Generally, a Finsler metric with
such a property is called an Einstein metric, i.e.,
\begin{equation}\label{k}
  \mathfrak{Ric}(x,y)=(n-1)K(x)
\end{equation}
for some function $K(x)$ on $M$. In particular, for a Randers
manifold $(M,F)$ with $\dim M\geq 3$, $F$ is an Einstein metric if
and only if there is a constant $K$ such that (\ref{k}) holds (see \cite{BR1}). The
following lemma is an important result on Einstein-Randers metrics.

\begin{lemma}[\cite{BR1}]\label{lem1}
Suppose $(M,F)$ is a Randers space with the navigation data
$(h,W)$. Then $(M,F)$ is Einstein with Ricci scalar
$\mathfrak{Ric}(x)=(n-1)K(x)$ if and only if there exists a constant
$\sigma$ satisfying the following conditions:
\begin{enumerate}
  \item $h$ is Einstein with Ricci scalar
  $(n-1)K(x)+\frac{1}{16}\sigma^2$, and
  \item $W$ is an infinitesimal homothety of $h$, i.e., $\mathfrak{L}_Wh=-\sigma
  h$.
\end{enumerate}
Furthermore, $\sigma$ must be zero whenever $h$ is not Ricci-flat.
\end{lemma}

It is well know that $K(x)$ is a
constant if $(M,F)$ is a homogeneous Einstein Finsler manifold. Here a Finsler manifold $(M,F)$ is called homogeneous if
its full group of isometries acts transitively on $M$. Based on
Lemma~\ref{lem1}, Deng-Hou obtained a characterization of
homogeneous Einstein-Randers metrics.

\begin{lemma}[\cite{DH1}]\label{lem2}
Let $G$ be a connected Lie group and $H$ a closed subgroup of $G$
such that $G/H$ is a reductive homogeneous space with a
decomposition ${\mathfrak g}={\mathfrak h}+{\mathfrak m}$. Suppose
$h$ is a $G$-invariant Riemannian metric on $G/H$ and $W\in
{\mathfrak m}$ is invariant under $H$ with $h(W,W)<1$. Let
$\widetilde{W}$ be the corresponding $G$-invariant vector field on
$G/H$ with $\widetilde{W}|_o=W$. Then the Randers metric $F$ with the
navigation data $(h,\widetilde{W})$ is Einstein with Ricci constant
$K$ if and only if $h$ is Einstein with Ricci constant $K$ and $W$
satisfies
\begin{equation}\label{asso}
  \langle [W,X]_{\mathfrak m}, Y\rangle+\langle X, [W,Y]_{\mathfrak
  m}\rangle=0, \quad \forall X,Y\in {\mathfrak m},
\end{equation}
where $\langle,\rangle$ is the restriction of $h$ on $T_o(G/H)\simeq
{\mathfrak m}$. In this case, $\widetilde{W}$ is a Killing vector
field with respect to the Riemannian metric $h$.
\end{lemma}

For simple, denote $H_i=A_4,{\mathfrak h}_j=A_1$, or $H_i=A_1,{\mathfrak h}_j=A_4$. By the equivalence of the adjoint representation and the isotropy representation of $H_i$ on
${\mathfrak h}_0\oplus{\mathfrak h}_j\oplus{\mathfrak m}_1\oplus{\mathfrak m}_2$, the
vector field
$$\widetilde{W}|_{gH}=d(\tau(g))|_o(W), \forall g\in {G}, W\in
{\mathfrak h}_0$$ is well-defined, and it is $G$-invariant (see
\cite{DH2}). For every metric given in Theorem~\ref{theorem1} and Theorem~\ref{theorem2}, one can
easily verifies the equation
$$ \langle [W,X]_{{\mathfrak h}_0\oplus{\mathfrak h}_j\oplus{\mathfrak m}_1\oplus{\mathfrak m}_2}, Y\rangle_{G/H_i}+\langle X, [W,Y]_{{\mathfrak h}_0\oplus{\mathfrak h}_j\oplus{\mathfrak m}_1\oplus{\mathfrak m}_2}\rangle_{G/H_j}=0$$
holds for any $W\in {\mathfrak h}_0$ and $X,Y\in {\mathfrak
h}_0\oplus{\mathfrak h}_j\oplus{\mathfrak m}_1\oplus{\mathfrak m}_2$, using the facts
that ${\mathfrak h}_0\subset {\mathfrak h}$ and that the metric is
$Ad(H)$-invariant. Then by Lemma~\ref{lem2}, the homogeneous metric
$$F(x,y)=\frac{\sqrt{[\langle W,y\rangle_{G/H_1}]^2+\langle y,y\rangle_{G/H_1}\lambda}}{\lambda}-\frac{\langle W,y\rangle_{G/H_1}}{\lambda}$$
is a $G$-invariant Einstein-Randers metric on $G/H_1$ when $W$ satisfies $\langle W,W\rangle_{G/H_1}<1$, and $F$ is Riemannian if and only if $W=0$. That is, we have the following theorem.
\begin{theorem}
There are at least four families of $E_6$-invariant non-Riemannian Einstein-Randers metrics on $E_6/A_4$, and two families of $E_6$-invariant non-Riemannian Einstein-Randers metrics on $E_6/A_1$.
\end{theorem}

\section{Acknowledgments}
This work is supported in part by NSFC (nos.11571182 and 11547122) and Nanhu Scholars Program for Young Scholars of XYNU.

\end{document}